\documentclass[12pt]{article}

\usepackage{GoetzArticle}
\usepackage{amssymb}
\usepackage{verbatim}
\usepackage{amssymb}
\usepackage{amsmath}
\usepackage[dvips]{color}
\usepackage{rotating}
\usepackage{bm}

\usepackage{rawfonts}
\input prepictex
\input pictexwd
\input postpictex

\newcommand{\commentout}[1]{}

\newcommand{\Z}{\mathbb{Z}}

\newcommand{\supp}{{\rm supp} \,}

\newcommand{\h}{{\cal H}}
\newcommand{\R}{\mathbb{R}}

\newcommand{\N}{\mathbb{N}}

\newcommand{\C}{\mathbb{C}}

\newcommand{\1}{\mathbf{1}}

\newcommand{\eqa}[1]{\begin{eqnarray}#1 \end{eqnarray}}

\newcommand{\Sc}{\mathcal S}
\newcommand{\RRh}{\R{\times}\widehat \R}

\newcommand{\Rh}{{\widehat\R}}

\newcommand{\dtrain}{{ \scriptstyle \mathrel\perp\joinrel\perp\joinrel\perp}}

\title{Measurement of time--varying Multiple--Input Multiple--Output Channels}
\author{
    G.~E. Pfander\footnotemark[3]
    }

\newcommand{\spacing}[1]{\renewcommand{\baselinestretch}{#1}\large\normalsize}
\spacing{1.333}

\begin{document}


\maketitle
\renewcommand{\thefootnote}{\fnsymbol{footnote}}

\footnotetext[3]{School of Engineering and Science,
    International University Bremen, 28759 Bremen, Germany.
    }

\renewcommand{\thefootnote}{\arabic{footnote}}
\pagenumbering{arabic} \maketitle

\begin{abstract} We derive a criterion on the measurability /
identifiability of Multiple--Input Multiple--Output (MIMO) channels
based on the size of the so-called spreading support of its
subchannels. Novel MIMO transmission techniques provide
high-capacity communication channels in time-varying environments
and exact knowledge of the transmission channel operator is of key
importance when trying to transmit information at a rate close to
channel capacity.

\end{abstract}







\section{Introduction}\label{sec:intro}


The recovery of information from a signal that has traveled through
a communications channel requires knowledge of --- or at least some
information on --- the transmission channel at hand.  In
applications such as mobile telephony, neither the location of the
subscriber nor the changing environment through which information is
transmitted is known {\it a-priori}. To combat this problem, a pilot
signal is send prior to information transmission with the hope that
the corresponding channel output supplies the receiver with the
measurements that are needed to invert the channel operator. The
inverse of the channel operator allows the receiver to recover the
information from the subsequently send information carrying signals.

In  Single--Input Single--Output (SISO) channels, the channel input
is considered to be a single variable function, which, after being
transmitted, is distorted by the unknown transmission channel
operator before arriving at the receiver (see \cite{Gol05,PNG03} and
references within). In \cite{PW06b}, the existence of pilot signals
which identify linear SISO channel operators was shown to depend on
the size of the spreading support of the channel operator. That is,
it was shown that a channel operator is identifiable by the channel
output corresponding to an appropriately chosen input signal if the
{\it a-priori} known spreading support has area (Jordan content)
less than one, while a channel operator cannot be identified by a
single input/output pair if the area of the spreading support is
larger than one (and nothing else is known of the channel operator).
Loosely speaking, the size of the spreading support of an operator
represents the amount of time--frequency dispersion that the channel
inflicts on the transmission signal. Too much time--frequency
dispersion cannot be resolved by a single channel output.
Fortunately, channel operators with spreading support area much
smaller than one, often called slowly time--varying or underspread
operators, are the norm in mobile communications. The results in
\cite{PW06b} described above were conjectured in the 1960s by
Kailath \cite{Kai62} and Bello \cite{Bel69}. See \cite{KP06} and
\cite{PW06b} for some historical background on the channel
identification problem for slowly time--varying channels and for
further applications of identification theorems for underspread
operators.

Multiple transmit and receive antenna methods have been developed to
obtain high capacity wireless channels (see
\cite{Gol05,Kue06,PNG03,Tel99} and references within). Methods which
achieve high capacities often rely on the precise knowledge of the
channel at the receiver and/or the transmitter (see \cite{Gol05}, pp
298).

In such MIMO channel setups, $N$ signals are transmitted by $N$
antennas simultaneously. On the receiver side, $M$ antennas record
channel output signals that represent the superposition of the $N$
input signals, each individually distorted depending on the path the
signal has travelled from its transmitting antenna to the receiving
antenna. Consequently, a linear MIMO  channel operator can be
modelled by a matrix of $N{\cdot}M$ SISO channel operators. It maps
a vector of $N$ transmission signals to $M$ channel output signals.

In this paper, we extend the SISO results from \cite{PW06b} to
linear MIMO channels. That is, we show that MIMO channel operators
permit identification by one vector of input signals if at each of
the receiving antennas the following condition holds: the sum of the
areas of the $N$ spreading supports of the subchannels leading to
the receiving antenna is less than one. Conversely, we show that if
the sum of the $N$ spreading areas of the subchannels leading to one
of the receiving antennas is larger than one, then identification is
not possible.

For simplicity, we assume throughout this paper that the $N{\cdot}M$
subchannels within a MIMO channel are independent of each other.
That is, information obtained on one of the $N{\cdot}M$ subchannels
does not carry any information on another subchannel in the MIMO
setup. The realistic assumption that the vicinity of the transmit
antennas and the vicinity of the receive antennas lead to a
dependent channel ensemble should allow for a relaxation of the
measurability criterion given here.

Modern methods in time--frequency analysis, such as those involving
Feichtinger's algebra and modulation spaces, have been used in
\cite{KP06,PW06,PW06b} to streamline the analysis of operators with
compactly supported spreading functions. Using these methods comes
at the price of necessitating non--standard terminology when
formulating results. Here, we bypass these methods in order to state
results in terms of the better known Hilbert--Schmidt operators and
tempered distributions. Further, the approach chosen here leads to a
generalization of the results in \cite{PW06b}  in the SISO case as
well.

Section~\ref{section:prelim} is devoted to preliminaries and
notation. We state our main result as Theorem~\ref{theorem:main} in
Section~\ref{section:mainresult}. The result is then proven in
Section~\ref{section:sufficient} and Section~\ref{section:necessary}


\section{Preliminaries and Notation}\label{section:prelim}

The space of complex valued Lebesgue integrable functions on
$d$-dimensional Euclidean space $\R^d$ is denoted by $L^1(\R^d)$.
The {\em Fourier transform} $\widehat f$ of $f\in L^1(\R^d)$ is the
continuous function
$$
\widehat f (\gamma) =\int f(x)\, e^{-2\pi i\gamma{\cdot} x}\, dx
\,,\quad \gamma \in \widehat \R^d,
$$
where $\widehat \R^d$ is the dual group of $\R^d$, which, aside of
notation, is identical to $\R^d$.

The space of square integrable functions $L^2(\R^d)$ consists of
those Lebesgue measurable functions which satisfy
$$
    \|f\|_{L^2} = \left( \int |f(x)|^2\, dx \right)^{\frac 1 2
    }<\infty.
$$
$L^2(\R^d)$ is a Hilbert space  with inner product
$$
\langle f, g \rangle =\int f(x)\overline{g(x)}\,dx,\quad f,g\in
L^2(\R^d).
$$
In case of vector valued functions $\bm{f}=(f_1,\ldots,f_N) \in
L^2(\R^d)^N$ we set accordingly
$$
    \|\bm{f}\|_{L^2}= \sqrt{\sum_{n=1}^N \|f_n\|^2_{L^2}}.
$$
For $f\in L^1(\R^d)\cap L^2(\R^d)$ we have $\|
\widehat f\|_{L^2}=\|
 f\|_{L^2}$. In fact, the Fourier transform on
$L^1(\R^d)\cap L^2(\R^d)$ extends to a unitary operator on the
Hilbert space $L^2(\R^d)$.

The set of {\em Schwartz class functions} $\Sc(\R)\subseteq L^2(\R)$
on $\R$ consists of all infinitely  differentiable functions which
satisfy
$$
    p_{k,l}(f) = \sup_{x\in\R} \big| x^l
    f^{(k)}(x)\big|<\infty,\quad k,l\in\N,
$$
where $f^{(k)}$ denotes the $k$-th derivative of $f$. $\Sc(\R)$ is a
Frech\'et space whose metric is defined using the seminorms
$p_{k,l}$, $k,l\in\N$. Hence, $f_n \to f$ in $\Sc(\R)$ if and only
if $p_{k,l}(f_n-f)\to 0$ for all $k,l\in\N$. The elements in the
dual space $\Sc'(\R)$ of bounded functionals on $\Sc(\R)$ are called
{\em tempered distributions}.

The usefulness of $\Sc(\R)$ and $\Sc'(\R)$ in harmonic analysis
stems in part from the fact that the Fourier transform defines a
bijective isomorphism on $\Sc(\R)$. Using duality, we can extend the
Fourier transform on $\Sc(\R)$ to the space $\Sc'(\R)$ of tempered
distributions. Since $\Sc'(\R)$ contains constant functions, {\em
Dirac's delta} $\delta:f\mapsto f(0)$, and {\em Shah distributions}
$\dtrain_a=\sum_{n\in \Z} \delta_{an}$, where $\delta_{na}=T_{na}
\delta$ and  $a>0$, it is justified to write $\widehat \dtrain_a =
\tfrac 1 a \dtrain_{\frac  1 a}$.

Similarly to the Fourier transform, the {\em time shift operator}
$T_t$, $t\in\R^d$, given by $T_tf(x)=f(x-t)$ and the {\em modulation
operator} $M_\omega$, $\omega \in \Rh^d$, $M_\omega f(x)=e^{2\pi i
\omega{\cdot} x } f(x)$ are unitary operators on $L^2(\R^d)$ and
bijective isomorphism on $\Sc(\R)$ and $\Sc'(\R)$ (equipped with the
weak-$\ast$ topology). Note that $M_\omega$ is also called {\em
frequency shift operator} since $\widehat{M_\omega f}=T_\omega
\widehat f$. Further, we refer to $\pi(\lambda)=\pi(t,\nu)=T_tM_\nu$
for $\lambda=(t,\nu)\in\R^d{\times}\Rh^d$ as {\em time--frequency
shift operator}.

The set $HS(L^2(\R))$ of {\em Hilbert--Schmidt operators} on
$L^2(\R)$ consists of those linear operators on $L^2(\R)$ which
satisfy
\begin{eqnarray}
  Hf(x)=\int \kappa_H(x,y)f(y)\, dy,\quad f\in\Sc(\R),
  \label{equation:HS}
\end{eqnarray}
for $\kappa_H\in L^2(\R^2)$ \cite{Die70,Gaa73}. In fact, the density
of $\Sc(\R)$ in $L^2(\R)$ together with $\langle Hf,g\rangle=\langle
\kappa_H, g\otimes \overline{f}\rangle$ implies that
(\ref{equation:HS}) extends to a bounded operator on $L^2(\R)$. Note
further, that $HS(L^2(\R))$ is  a Hilbert space with inner product
$\langle H_1,H_2\rangle_{HS}=\langle
\kappa_{H_1},\kappa_{H_2}\rangle$ and corresponding norm.
Hilbert--Schmidt operators are compact operators on $L^2(\R)$. Note
that  some Hilbert--Schmidt operators can be extended to act on
larger subsets of $\Sc'(\R)$ than $L^2(\R)$, a fact that will use
later in this paper.

Every Hilbert--Schmidt operator can be expressed as a superposition
of time and frequency shift operators. In fact, for $H$ with
$\kappa_H \in L^1(\R^2)\cap L^2(\R^2)$, we set
$$
    \eta_H(t,\nu)=\int \kappa_H(x,x-t) e^{-2\pi i \nu x }\, dx,\quad
    a.e.\ \nu\in\widehat R.
$$
It is easy to see that in this case
\begin{eqnarray}
  \|\eta_H\|_{L^2}=\|\kappa_H\|_{L^2}=\|H\|_{HS},
  \label{equation:eta=HSnorm}
\end{eqnarray} implying that the {\em spreading function}
$\eta_H\in L^2(\R{\times}\Rh)$ can be defined for any
Hilbert--Schmidt operator $H$, and thereby extending
(\ref{equation:eta=HSnorm}) to all Hilbert--Schmidt
operators\footnote{The spreading function of an Hilbert--Schmidt
operator, or, more general, of a pseudodifferential operator, is the
symplectic Fouriertransform of the operators Kohn--Nirenberg symbol.
Consequently, the theory of pseudodifferential operators with
compactly supported spreading functions coincides with the theory of
pseudodifferential operators with bandlimited Kohn--Nirenberg
symbols.}. As mentioned above, we have
\begin{eqnarray}
  H =\iint \eta_H(t,\nu) T_t M_\nu \, d\nu dt=\int
  \eta_H(\lambda)\pi(\lambda)\,d\lambda
  \label{equation:weaketa}
\end{eqnarray}
where the operator valued integral in (\ref{equation:weaketa}) is
understood weakly, that is, $H$ is defined via
\begin{eqnarray}
    \langle Hf,g\rangle=\iint\eta_H(t,\nu) \int e^{2\pi i \nu
(x-t)}f(x-t)\overline{g(x)}dx \, dt d\nu=\langle \eta_H, V_f g
\rangle, \label{equation:ExtendHS}
\end{eqnarray}
where the {\em short--time Fourier transform} $V_fg$ of $g\in
L^2(\R)$ with respect to $f\in L^2(\R)$ is given by
$$V_f g(t,\nu)=\int g(x) e^{- 2\pi i \nu
(x-t)}\overline{f(x-t)}\,dx$$ and satisfies $V_f g\in
L^2(\R{\times}\widehat \R)$  \cite{Gro01}.

To avoid double indices, we shall write at times $\eta(H)$ in place
of $\eta_H$ and, similarly, $\kappa(H)$ in place of $\kappa_H$.

We denote by $HS(L^2(\R))^{M{\times}N}$ the space of $N$-input,
$M$-output MIMO channels whose $N{\cdot}M$ subchannels are
Hilbert--Schmidt operators on $L^2(\R)$ \cite{GP06}. The operator
space $HS(L^2(\R))^{M{\times}N}$ is equipped with norm
$$
    \big\|\bm H\big\|_{HS}
        =\sqrt{\sum_{m=1}^M\sum_{n=1}^N \big\|H_{mn}\big\|_{HS}^2},\quad
    \bm H=\left(%
\begin{smallmatrix}
  H_{11} & \cdots & H_{1N}  \\
  \vdots &  & \vdots \\
  H_{M1} & \cdots & H_{MN} \\
\end{smallmatrix}%
\right) \in HS(L^2(\R))^{M{\times}N}.\ \footnote{It is easy to see
that $HS(L^2(\R))^{N{\times}N}=HS(L^2(\R)^N)$.}
$$
Further, the spreading function $\bm \eta_{\bm H} =\bm \eta(\bm H)$
of $\bm H=\left(%
\begin{smallmatrix}
  H_{11} & \cdots & H_{1N}  \\
  \vdots &  & \vdots \\
  H_{M1} & \cdots & H_{MN} \\
\end{smallmatrix}%
\right)\in  HS(L^2(\R))^{M{\times}N}$ and the {\em spreading
support} of $\bm H$ are defined componentwise, that is, we have
$$\bm \eta(\bm H)= \left(%
\begin{smallmatrix}
  \eta(H_{11}) & \cdots & \eta(H_{1N})  \\
  \vdots &  & \vdots \\[.2cm]
  \eta(H_{M1}) & \cdots & \eta(H_{MN}) \\
\end{smallmatrix}%
\right)\in L^2(\RRh)^{M\times N},$$ and $$ \supp\, \bm \eta(\bm
H)=\left(\begin{smallmatrix}
  \supp\, \eta(H_{11}) & \cdots & \supp\, \eta(H_{1N})  \\
  \vdots &  & \vdots \\[.2cm]
  \supp\, \eta(H_{M1}) & \cdots & \supp\, \eta(H_{MN}) \\
\end{smallmatrix}%
\right)\subseteq (\RRh)^{M\times N}.$$


Our identifiability result for MIMO channels considers operator
classes of the form
$$
\h_{\bm S}=\Big\{\bm H\in HS(L^2(\R))^{M{\times}N}:\ \supp\,\bm
\eta(\bm H)\subseteq \bm S \Big\},\quad \bm S\subseteq
(\RRh)^{M\times N}.
$$
To avoid pathological cases, we shall only consider $\h_{\bm S}$
where $\bm S$ is the cartesian products of so called Jordan domains.

\begin{definition}
A {\em Jordan domain} $M\subseteq \RRh$ is a bounded set whose
boundary is a Lebesgue zero set.
\end{definition}

Clearly, our restriction to Jordan domains is not relevant to
applications such as those in communications engineering. The
following useful characterization of Jordan domains is well known.
It is discussed in detail in \cite{LPW05}.

\begin{lemma}
If $M$ is a Jordan domain, then its Lebesgue measure $\mu(M)$
satisfies
\begin{eqnarray}
    \mu(M)&=& \sup\{ \mu (U): U\subseteq M\text{ and } U\in \mathcal
        U_{KL} \text{ for some } K,L\in\N,\ L \text{ prime } \}\notag \\
        &=& \inf\{ \mu(U): U\supseteq M\text{ and } U\in \mathcal U_{KL}
        \text{ for some } K,L\in\N, \ L \text{ prime } \}.\notag
\end{eqnarray}
where for $K,L\in\N$ we  set $R_{KL}=[0,\tfrac 1 K
]\times[0,\tfrac{K}{L}]$ and
\begin{eqnarray*}
  \mathcal{U}_{KL}
        =   \left\{
                \bigcup_{j=1}^J \left(R_{KL}+(\tfrac {m_j} K  , \tfrac {n_j K} {L}   )
                \right):
                \ m_j,n_j\in\Z, J\in\N
            \right\}. \label{equation:U-KL}
\end{eqnarray*}
\end{lemma}

\section{Statement of Results}\label{section:mainresult}

The domain of Hilbert--Schmidt operators with compactly supported
spreading function can be extended to include classes of tempered
distributions (see Theorem 4.2 in \cite{PW06}).  For example, using
(\ref{equation:ExtendHS}), it is easy to see that any
Hilbert--Schmidt operator with compactly supported spreading
function maps $\dtrain_a$, $a\in \R^+$, to a function in $L^2(\R)$.
In fact, a simple computation in \cite{KP06} shows that for
$S=[-\frac 1 2, \frac 1 2]{\times}[-\frac 1 2, \frac 1 2]\subseteq
\RRh $ we have
$$\|H\dtrain_1\|_{L^2(\R)}=\|H\|_{HS}, \quad H\in\h_S. $$

\begin{definition}
An operator class $\h \subseteq HS(L^2(\R))^{M{\times}N}$  is
identifiable if there exists $\bm f\in S'(\R)^N$ and positive $A,B$,
with $$A\,\|\bm H\|_{HS}\leq \|\bm Hf\|_{L^2}\leq B\, \|\bm
H\|_{HS}\quad \text{for } \bm H\in \h.$$
\end{definition}

In short, an operator class $\h$ is identifiable if there is $\bm f$
with the property that the induced map
$$
\Phi_{\bm f}:\h\longrightarrow L^2(\R)^N,\
    \bm H\mapsto \bm H \bm f
$$
is bounded and stable, that is, bounded above and below.

\begin{theorem}\label{theorem:main}
Let $\bm S=(S_{mn})\subseteq (\RRh)^{M\times N}$ be the cartesian
product of Jordan domains in $\RRh$ and let
$$
\h_{\bm S}=\Big\{\bm H\in HS(L^2(\R))^{M{\times}N}:\
\supp\,\bm\eta(\bm H)\subseteq \bm S \Big\}.
$$
\begin{enumerate}
  \item
  If $\displaystyle\sum_{n=1}^N \mu(S_{mn})<1$
  for all $m\in\{1,\ldots,M\}$, then
  $\h_{\bm S}$ is identifiable.

  \item
  If $\displaystyle\sum_{n=1}^N \mu(S_{mn})>1$
  for some $m\in\{1,\ldots,M\}$, then
  $\h_{\bm S}$ is not identifiable.
\end{enumerate}

\end{theorem}

\section{Proof of Theorem 3.2, part {\it
1}}\label{section:sufficient}

Theorem~\ref{theorem:sufficientidentifiabilitywalnut} reduces
Theorem \ref{theorem:main}, part {\it 1}, for SISO channels
($M=N=1$) to a question on the linear independence of columns of the
following matrices: for any $L$--periodic sequence $c=\{c_k\}_{k\in
Z}$ we set $\bm A(c) = [\bm A_0(c)\,\,\,\bm
A_1(c)\,\,\,\cdots\,\,\,\bm A_{K-1}(c)]\in\C^{KL{\times}L}$ with
$\displaystyle \bm A_k(c) = (c_{p+k}\,e^{2\pi
iq(p+k)/L})_{p,q=0}^{L-1}\in\C^{L{\times}L}$.

\begin{theorem}\label{theorem:sufficientidentifiabilitywalnut}
Let $c=\{c_k\}_{k\in Z}$ be a sequence with period $L$ and
$\displaystyle f = \sum_k c_k\,\delta_{\frac k K}\in \Sc'(\R).$
Further, set
$$U=
                \bigcup_{j=1}^J \left(R_{KL}+(\tfrac {m_j} K  , \tfrac {n_j K} {L}   )
                \right),\quad
                \ m_j,n_j\in\Z, J\in\N,$$ where $R_{KL}=[0,\frac 1 K]\times[0,\frac K L]$.

Then $f$ identifies $\mathcal H_U$ if and only if the columns in
$\bm A(c)$ with column indices in $\{m_j L + n_j\}_j$ are linearly
independent.
\end{theorem}

Clearly, this result is only applicable if the cardinality $|J|$ of
$J$ satisfies $|J|\leq L$ since $\bm A(c)$ has at most $L$ linear
independent columns. This requirement is equivalent to $\mu(U)\leq
|J| \frac 1 K \frac K L \leq 1$.

If $L$ is prime, then $|J|\leq L$   is also sufficient for the
existence of an identifier for a SISO channel\cite{LPW05}:

\begin{theorem}\label{theorem:HaarExists}
If $L$ is prime then there exists $c\in\C^L$ such that any set of
$L$ columns of $\bm A(c)$ is linearly independent.
\end{theorem}


{\it Proof of Theorem \ref{theorem:main}, Part 1.}

We choose $\bm S=(S_{mn})\subseteq (\RRh)^{M{\times}N}$ which
satisfies $\displaystyle \sum_{n=1}^N \mu(S_{mn}) <1$. Since all
$S_{mn}$ are assumed to be Jordan domains, there exists $K, L\in\N$,
$L$ prime, so that for each $S_{mn}$ exists $U_{mn} \in \mathcal
U_{KL}$ with $S_{mn}\subseteq U_{mn}$ and $\displaystyle
\sum_{n=1}^N \mu(U_{mn}) <1$ for $m=1,\ldots,M$.

Clearly, $\h_{\bm S}\subseteq \h_{\bm U}$with $\bm
U=(U_{mn})\subseteq (\RRh)^{M\times N}$ implies that the
identifiability of $\h_{\bm S}$ follows from the identifiability of
$\h_{\bm U}$ which we shall prove now.

All $U_{mn}$ are bounded, hence, we can choose $W>0$ so that
$$U_{mn}\subseteq B^\infty_W(0)=\Big\{\|(t,\nu)\|_\infty=\max \{
|t|,|\nu|\} \leq W\Big\}\quad \text{for } m=1,\ldots, M,\
n=1,\ldots,N.$$ For $L$ and $K$ chosen above,
Theorem~\ref{theorem:HaarExists} allows us to choose an
$L$--periodic sequence $c$ so that any set of $L$ columns from $\bm
A(c)$ is linearly independent. We set
$$
    f_n=\pi(0,(n{-}1)2W)\sum_{k\in\Z} c_{k\text{mod}L}
    \delta_{\frac k K}\quad \text{for } n=1,\ldots, N\,,
$$
and claim that $\bm f=(f_1,\ldots,f_N)^T$ identifies $\h_{\bm U}$.

To see this, note that the choice of $W$ implies that
$T_{(0,(n{-}1)2W)}U_{mn}\cap T_{(0,(n'-1)2W)}U_{mn'}=\emptyset$ for
all $n\neq n'$ and $m=1,\ldots , M$. For $U_m= \bigcup_{n=1}^N
T_{(0,(n{-}1)2W)}U_{mn}$, $m=1,\ldots, M$, we have
$\mu(U_m)=\displaystyle\sum_{n=1}^N \mu(U_{mn})<1$, and, by
Theorem~\ref{theorem:sufficientidentifiabilitywalnut}, $f_1$
identifies $\h_{U_m}\subseteq HS(\R)$ for $m=1,\ldots,M$, that is,
there exists $A,B>0$ such that for all $H\in \h_{U_m}$,
$m=1,\ldots,M$ we have
\begin{eqnarray}
    A\|H\|_{HS}= A\|\eta_H\|_{L^2}\leq \|Hf_1\|_{L^2}\leq
B\|H\|_{HS}.\label{equation:identificationUm}
\end{eqnarray}
For $\bm H\in\h_{\bm U}$ we set $\bm g=(g_1,\ldots,g_M)=\bm H\bm f$
and compute for $m=1,\ldots, M$,
\begin{eqnarray}
  g_m   &=& \sum_{n=1}^N H_{mn}\ f_n =\sum_{n=1}^N H_{mn}\circ\pi(0,(n{-}1)2W)
            \ f_1\notag\\
        &=& \sum_{n=1}^N \int \eta\big(H_{mn}\circ
        \pi(0,(n{-}1)2W)\big)(\lambda)\ \pi(\lambda)
            f_1 \,d\lambda \notag\\
        &=&  \int \left(
                        \sum_{n=1}^N \eta\big(H_{mn}\big)(\lambda-(0,(n{-}1)2W))
                    \right)\ \pi(\lambda)
            f_1 \,d\lambda \notag.
\end{eqnarray}
Since $\supp T_{(0,(n{-}1)2W)}\eta\big(H_{mn}\big)\subseteq
T_{(0,(n{-}1)2W)}U_{mn}\subseteq U_m$ and
$$\mu\Big(\ \supp T_{(0,(n{-}1)2W)}\eta\big(H_{mn}\big) \cap
\supp T_{(0,(n'{-}1)2W)}\eta\big(H_{mn'}\big)\ \Big)=0$$ for all
$n\neq n'$ and all $m=1,\ldots, M$, we can apply
(\ref{equation:identificationUm}) to obtain
\begin{eqnarray}
  \|g_m\|^2_{L^2}&\geq& A^2 \left\|\sum_{n=1}^N T_{(0,(n{-}1)2W)}
  \eta\big(H_{mn}\big)
  \right\|^2_{L^2}=A^2 \sum_{n=1}^N\left\|\eta\big(H_{mn}\big)
  \right\|^2_{L^2}\notag
\end{eqnarray}
and
\begin{eqnarray}
  \|\bm g\|^2_{L^2}&=& \sum_{m=1}^M \|g_m\|^2_{L^2}
  \geq  A^2  \sum_{m=1}^M \sum_{n=1}^N\left\|\eta\big(H_{mn}\big)
  \right\|^2_{L^2}
  = A^2  \left\|\bm H
  \right\|^2_{HS}\notag
\end{eqnarray}
The upper bound involving $B$ follows in the same manner. \hfill
$\square$\par

\section{Proof of Theorem 3.2, part {\it
2}}\label{section:necessary}

We shall now show that the condition $\displaystyle\sum_{n=1}^N
\mu(S_{mn})\leq 1$, $m=1\ldots M$, is necessary for the
identifiability of $\h_{\bm S}$, $\bm S= (S_{mn})$.

Without loss of generality, we  assume  a Multiple--Input
Single--Output (MISO) scenario, that is, we consider $M=m=1$ and
write $S_n=S_{1n}$ and $H_n=H_{1n}$. In fact, if there there exists
$\bm S$ in the MIMO case with $\displaystyle\sum_{n=1}^N
\mu(S_{m_0n})>1$ for $m_0\in\{1,\ldots,M\}$ and $\h_{\bm S}$
identifiable, then defining $\bm S'$ by $S'_n=S_{m_0n}$ would lead
to a contradiction of Theorem~\ref{theorem:main}, part {\it 2}, in
the MISO case.

The proof of Theorem~\ref{theorem:main}, part {\it 2}, is organized
as the corresponding proof in \cite{PW06}. The crux is to show that
operators  in the class $\h_{\bm S}$ with $\displaystyle\sum_{n=1}^N
\mu(S_{n})>1$ carry to many, in time and frequency tightly packed,
degrees of freedom, that is, too much information to be embedded in
a stable manner in a single output signal.

To see this, we shall fix $\bm S$ with $\displaystyle\sum_{n=1}^N
\mu(S_{n})>1$. For this $\bm S$, we construct a bounded and stable
synthesis (information embedding) map $E:l_0(\Z^2)\longrightarrow
\h_S $ where $l_0(\Z^2)$ is equipped with the $l^2(\Z^2)$-norm, and
a bounded and stable analysis (information recovery) operator
$C:L^2(\R)\longrightarrow l^2(\Z^2)$ with the property that {\bf
all} compositions
$$
    C\circ \Phi_{\bm f} \circ E :\ l_0(\Z^2) \longrightarrow l^2(\Z^2)
    ,\quad \bm f\in \Sc'(\R)^N,$$
are not stable. The stability of $E$ and $C$ implies that the
boxed-in operators $\Phi_{\bm f}:\h_S\longrightarrow L^2(\R)$, $\bm
f\in \Sc'(\R)^N$, must not be stable, showing that $\h_{\bm S}$ is
not identifiable if $\displaystyle\sum_{n=1}^N \mu(S_{n})>1$.

Before proving Theorem~\ref{theorem:main}, part {\it 2}, we state
three lemmas, some of whose proofs can be found in \cite{PW06}.
Lemma~\ref{lemma:time-frequency-spreading-shift} concerns the
conjugation of Hilbert--Schmidt operators by time--frequency shifts.
In Lemma~\ref{lem:elementarybound} we construct a prototype operator
which is later used to construct a Riesz bases for its closed linear
span in $\h_S$, that is, a family of Hilbert--Schmidt operators
$\{H_{k,l}\}_{k,l\in\Z}$ for which the map
$$
\begin{array}{crcl}
 E:&l^2(\Z^2)&\rightarrow &HS(L^2(\R))\\
   &\{c_{kl} \}_{k,l\in\Z}& \mapsto   &
  \sum_{k,l\in\Z}c_{k,l} H_{k,l}
\end{array}
$$
is well defined, bounded, and stable. Lemma~\ref{lem:matrixinverse}
generalizes the fact that $m\times n$ matrices with $m<n$ have a
nontrivial kernel and, therefore, are not stable, to operators
acting on $l^2(\Z^2)$. In fact, the bi-infinite matrices
$M=(m_{j',j})_{j',j\in\Z^2}$ considered in
Lemma~\ref{lem:matrixinverse} are not dominated by its diagonal
$m_{j,j}$ --- which would correspond to square matrices --- but by a
slanted diagonal $m_{j,\lambda j}$, $j\in\Z^2$, with $\lambda>1$.
%
\begin{lemma} \label{lemma:time-frequency-spreading-shift}
For $P\in HS(\R)$ with spreading function $\eta_P \in L^2(\RRh)$ set
$\widetilde{P}= M_{\omega}T_{p - r} P T_{r}M_{\xi-\omega}\in HS$.
Then $\eta_{\widetilde{P}}=e^{2\pi i \omega
p}M_{(\omega,r)}\,T_{(p,\xi)}\,\eta_P$ and $\widetilde{P}\in
HS(\R)$.
\end{lemma}

\begin{lemma}\label{lem:elementarybound}
Fix $\lambda>1$ with $1<\lambda^4<\mu(S)$ and choose even functions
$\eta_1,\eta_2\in \mathcal{S}(\R)$ with values in $[0,1]$ and
$$
  \eta_1(t)=
  \begin{cases}
    1 & \text{for}\ |t|   \leq \tfrac{1}{2\lambda K} \\
    0 & \text{for}\ |t|   \geq \tfrac{1}{2K}
  \end{cases}\quad \text{and}
  \ \quad
  \eta_2(\nu)=
  \begin{cases}
    1 & \text{for}\ |\nu|   \leq \tfrac{K}{2\lambda L}\\
    0 & \text{for}\ |\nu|   \geq \tfrac{K}{2L}
  \end{cases}\ .
$$
The operator  $P\in \h_{R_{KL}}$ defined by $\eta_P=\eta_1\otimes
\eta_2$ has the properties:

\noindent a) The operator family
\begin{eqnarray}
   \left\{
            M_{\lambda  K k}\,
            T_{\frac 1 K m -\frac  {\lambda L}{K}l }\,
            P\,
            T_{\frac {\lambda L}{ K}l} \,
            M_{\frac{ K}{L}n-\lambda K k}\,\right\}_{k,l,m,n\in\Z}
\label{equation:OperatorFamily}
\end{eqnarray}
is a Riesz basis for its closed linear span in the Hilbert space
of Hilbert--Schmidt operators $HS(\R)$.

\noindent b) For $f\in \Sc'(\R)$, there exists $C_f,L_f\in\N$ and
$d_1,d_2:\R\rightarrow\R_0^+$ which decay rapidly at infinity with
$$ |PT_yM_\omega f(x)|\leq C_f \,d_1(x)(1+\|(y,\omega)\|_\infty)^{L_f},\quad x\in\R,$$
and
$$
 |\widehat{PT_yM_\omega f}(\xi)|\leq C_f\,d_2(\xi)(1+\|(y,\omega)\|_\infty)^{L_f},\quad
 \xi\in\widehat{\R}.
 $$
\end{lemma}

\begin{proof} {\rm (a)} See \cite{KP06}.

{\rm (b)}  For $f\in \Sc(\R)$, we compute
$$
    Pf(x)
        = \iint\eta_1(t)\eta_2(\nu)e^{2\pi i \nu (x-t)} f(x-t)\, d\nu\,dt
        = \int \eta_1(t) \check{\eta_2}(x-t) f(x-t)\,dt
        = \eta_1\ast (\check{\eta_2}f),
$$
and, therefore, $ \widehat{Pf}(\xi)=\widehat{\eta_1}(\xi)\cdot
\eta_2{\ast}\widehat{f}(\xi)$. The rapid decay and smoothness of
$\widehat \eta_1$ together with the fact that $\supp \eta_2$ compact
and $\eta_2$ smooth implies that $\widehat{Pf}$ and, therefore, $Pf$
is well defined for $f\in \Sc'(\R)$. In fact, we can conclude that
$\widehat {Pf}$, and, therefore, $Pf\in \Sc(\R)$ for $f\in
\Sc'(\R)$.

Further, we obtain for $f\in \Sc(\R)$ and $\xi\in\Rh$ that
\begin{eqnarray*}
\left| (PT_{-y}M_{-\omega}f)\widehat{\ }\,(\xi)\right|
    &=&\left|\widehat{\eta_1}(\xi)\,
        \int \eta_2(\xi - \nu)\, M_{-y}T_{\omega}\widehat{f}(\nu)d\nu\right|
    =|\widehat{\eta_1}(\xi)|\,\left|\langle M_{-y}T_{\omega}\widehat{f},\, T_\xi \eta_2\rangle \right|
        \\
    &=&|\widehat{\eta_1}(\xi)|\,\left|\langle \widehat{f},\, M_y T_{\xi-\omega} \eta_2 \rangle \right|
    =|\widehat{\eta_1}(\xi)|
        \left| V_{\eta_2} \widehat f (\xi-\omega,y) \right|.
\end{eqnarray*}
The weak-$\ast$ density of $\Sc(\R)$ in $\Sc'(\R)$ extends the
equality above to $f\in \Sc'(\R)$. Theorem 11.2.3 in \cite{Gro01}
provides us now  with $C_f',L_f'\in\N$ and
\begin{eqnarray*}
\left| (PT_{-y}M_{-\omega}f)\widehat{\ }\,(\xi)\right|
    &=&|\widehat{\eta_1}(\xi)|
        \left| V_{\eta_2} \widehat f (\xi-\omega,y) \right|
    \leq C_f'|\widehat{\eta_1}(\xi)| (1+|y|+|\xi-\omega|)^{L_f'}\\
    &\leq&C_f'|\widehat{\eta_1}(\xi)|\,(1+|y|+|\xi|+|\omega|)^{L_f'}\\
    &\leq&C_f'|\widehat{\eta_1}(\xi)|\,(1+|(\xi)|)^{L_f'} (1+|y|+|\omega|)^{L_f'}
    \leq d_2(\xi) (1+\|(y,\omega)\|_\infty)^{L_f'},
\end{eqnarray*}
where $d_2=C_f'\, 2^{L_f'}\,
|\widehat{\eta_1}(\xi)|\,(1+|\xi|)^{L_f'}$ is rapidly decaying.

Similarly, we conclude that for $f\in \Sc(\R)$ and $x\in\R$ we have
\begin{eqnarray*}
\left| PT_{-y}M_{-\omega}f(x)\right|
    &=&\left| \int \eta_1(s-x) \check{\eta_2}(s)\, T_{-y}M_{-\omega}f(s)\,ds   \right|
        \\
    &=&\left|\langle M_{\omega}T_{y}(\check{\eta_2}T_x\eta_1) ,\, f\rangle \right|
       =
        \left| V_{\check{\eta_2}T_x\eta_1} f (y,\omega) \right|
\end{eqnarray*}
Within the proof of Theorem 11.2.3 in \cite{Gro01},  the existence
of $C_f,L_f\in\N$ are given with $C_f\geq C_f'$, $L_f\geq L_f'$, and
\begin{eqnarray}
\left| PT_{-y}M_{-\omega}f(x)\right|
    =  \left| V_{\check{\eta_2}T_x\eta_1} f (y,\omega) \right| \leq C_f
        \max_{m,n\leq L_f}\ \sup_{t\in\R}\ |t^n \frac{\partial^n}{\partial t^n
}\check{\eta_2}T_x\eta_1(t)|\ (1+\|(y,\omega)\|_\infty)^{L_f}.\notag
\end{eqnarray}
 Note that since
$\check{\eta_2},\eta_1\in\Sc(\R)$, each $\sup_{t\in\R}\ |t^n
\frac{\partial^n}{\partial t^n}\check{\eta_2}T_x\eta_1(t)|$,
$m,n\leq L_f$, decays faster than any polynomial. This implies that
also $d_1(x)=\max_{m,n\leq L_f}\ \sup_{t\in\R}\ |t^n
\frac{\partial^n}{\partial t^n }\check{\eta_2}T_x\eta_1(t)|$ also
decays faster than any polynomial.
\end{proof}

%
\begin{lemma}\label{lem:matrixinverse}
Given $M=(m_{j',j}):l^{2}(\Z^2)\rightarrow l^{2}(\Z^2)$. If there
exists a polynomial $p$ of degree $L\in\N$ and a monotonically
decreasing function $w:\R_0^+\rightarrow\R_0^+$ with
$w(x)=o\left(x^{-(L{+}2)}\right)$ satisfying
$$\displaystyle
  |m_{j',j}|< w(\|\lambda j'-j\|_\infty)\ p(\|j\|_\infty),\quad \|\lambda j'-j\|_\infty > K_0
$$
for some constants $\lambda >1$ and $K_0>0$, then $M$ is not stable.
\end{lemma}
The proof of Lemma~\ref{lem:matrixinverse} is included in the
appendix.


Now all pieces are in place to prove necessity of the condition
$\displaystyle\sum_{n=1}^N \mu(S_{n})\leq 1$ for the identifiability
of $\h_{\bm S}$, $\bm S=(S_{mn})$.

{\it Proof of Theorem~\ref{theorem:main}, part 2.}

Fix $\bm S=(S_n)$ with $\sum_{n=1}^N \mu(S_{n})>1$. Without
restriction of generality, we shall assume that $S_n\in
\mathcal{U}_{KL}$ for some $K,L\in\N$ and all $n=1,\ldots,N$, and
that $S_n\cap S_{n'}=\emptyset$ for $n\neq n'$. Hence, there exists
$\mathcal J=\{0,1,2,...,J-1\} \subseteq
            \N$ so that $S=\bigcup_{n=1}^N S_n= \bigcup_{j\in{\mathcal J}} \left(R_{KL}+(\tfrac
            {m_j} K , \tfrac {n_j K} {L}) \right)$, $(m_j,n_j)\neq
            (m_{j'},n_{j'})$ for $j\neq j'$.  We have $\mu(R_{KL})= \tfrac 1
            {L}$, and, since $  \mu(S)=\sum_{n=1}^N \mu(S_{n})>1$,
            we have $J > L$.

Fix $\bm f =(f_1,\ldots,f_n)\in \Sc'(\R)^N$. Choose $\lambda$,
$\eta_1$, $\eta_2$, $P$ , $C_{\bm f}=\max_n C_{f_n}$, $L_{\bm
f}=\max_n L_{f_n}$, and  to $C_{\bm f}$ and $L_{\bm f}$
corresponding $d_1$ and $d_2$ according to
Lemma~\ref{lem:elementarybound}. For $n=1,\ldots, N$ define
$$J_n=\left\{j\in\{0,\ldots,J-1\}:\ R_{KL}+(\tfrac
            {m_j} K , \tfrac {n_j K} {L})\subseteq S_n\right\}.$$

The synthesis operator $E:l_0(\Z^2)\rightarrow \h_S$ mentioned above
is given by
\begin{eqnarray}
  E:\sigma_{k,l''}= \sigma_{k,lJ+j} \mapsto \sum_{k,l\in Z}\sum_{j=0}^{J-1}\ \sigma_{k,lJ{+}j} \
            \iota(j)
            M_{ \lambda  K k}\,
            T_{\frac{1}{K}m_j + \frac  {\lambda L}{K}l }\,
            P\,
            T_{- \frac  {\lambda L}{K}l }\,
            M_{\frac{ K}{L}n_j-\lambda K k},\notag
\end{eqnarray}
where $$\iota(j):HS(\R)\longrightarrow HS(\R)^N,\quad H\mapsto
H\cdot (\1_{J_1}(j),\ldots,\1_{J_M}(j))=\stackrel{n^{th}\text{
position if $j\in J_n$}}{(0,\ldots 0, H,0,\ldots,0)}.$$ Since
\begin{eqnarray}
  \left\{
            M_{\lambda  K k}\,
            T_{\frac 1 K m -\frac  {\lambda L}{K}l }\,
            P\,
            T_{\frac {\lambda L}{ K}l} \,
            M_{\frac{ K}{L}n-\lambda K k}\right\}_{k,l,m,n\in\Z}
\notag
\end{eqnarray}
is a Riesz basis for its closed linear span in $\h_S\subseteq
HS(\R)$, we have that
\begin{eqnarray}
    \left\{
            \iota(j)M_{\lambda  K k}\,
            T_{\frac{1}{K}m_j +\frac  {\lambda L}{K}l }\,
            P\,
            T_{- \frac  {\lambda L}{K}l }\,
            M_{\frac{ K}{L}n_j-\lambda K k}\right\}_{k,l\in\Z,j\in\mathcal J}
\notag
\end{eqnarray}
is a Riesz basis for its closed linear span in $HS(\R)^N$. We
conclude that $E$ is bounded and stable.

To construct a stable analysis operator $C$, we choose the Gaussian
$g_0:\R\rightarrow\R^+,\quad x\mapsto e^{-\pi x^2}$, and note that
Lyubarski \cite{Lyu92} and Seip and Wallsten \cite{Sei92b,SW92} have
shown that $ \left\{M_{ka'}T_{lb'} g_0\right\}$ is a frame whenever
$a'b'<1$.\footnote{For background on frame theory see
\cite{Chr03,Gro01}.} Since $\lambda^2 K\frac{\lambda^2 L}{KJ}
=\lambda^4 \frac{L}{J} =\frac{\lambda^4}{ \mu(S)}<1,$ this implies
that the analysis map given by
$$
  C:L^2(\R)  \rightarrow  \l^2(\Z^2), \quad
         f\mapsto
         \left\{
           \langle f, M_{\lambda^2 K\,k}
            T_{\frac{\lambda^2  L}{KJ}\,l } g_0\rangle
         \right\}_{k,l}
$$
is bounded  and stable.

For simplicity of notation, set $\alpha =K$ and $ \beta=\frac L {KJ}
$. Let us now  consider the composition
$$
\begin{array}{ccccccc}
  l_0(\Z^2)      & \stackrel{E}{\rightarrow} & \h_{\bm S}
    & \stackrel{\Phi_{ \bm f} }{\rightarrow} &
  L^2(\R)
        & \stackrel{C}{\rightarrow} & l^2(\Z^2)
  \\
  \{{\sigma}_{k,l''}\} & \mapsto         & E{\{{\sigma}_{k,l''}\}} &  \mapsto   & E{\{{\sigma}_{k,l''}\}} \,\bm f
        & \mapsto         &
        \left\{\, \langle\, E\{\sigma_{k,l''}\}\,\bm f ,\ M_{\lambda^2\alpha k'}T_{\lambda^2\beta l'} \,
         g_0 \, \rangle \, \right\}_{k',l'}\, .
\end{array}
$$
We set $f_j=f_n$ whenever $j\in J_n$ and note that the bi--infinite
matrix
$$
M=\Big(m_{k',l',k,l''}\Big)=\Big(m_{k',l', \,k \,,lJ+j}\Big)=
    \Big(   \langle\, \
            M_{\lambda \alpha k}\,
            T_{\frac{m_j}{\alpha} + \lambda   \beta l J}\,
            P\,
            T_{- \lambda \beta l J  }\,
            M_{\frac{  n_j}{\beta J}-\lambda \alpha k}\,f_j\, ,\
                M_{\lambda^2\alpha k'}T_{\lambda^2 \beta l'} \, g_0\,
            \rangle
    \Big) \, ,
$$
$l''=lJ+j$, represents the operator $C\circ \Phi_{\bm f} \circ E$
with respect to the canonical basis of $l^2(\Z^2)$, since
\begin{eqnarray}
  \Big(C\circ \Phi_{\bm f} \circ E \ \{\sigma_{k,lJ+j}\}\Big)_{k',l'}
    &=&
    \langle\,
    \sum_{k,l}\sum_{j=0}^{J-1}
      \sigma_{k,lJ+j}\             M_{\lambda \alpha k}\,
            T_{\frac{m_j}{\alpha} + \lambda   \beta l J}\,
            P\,
            T_{-\lambda \beta l J  }\,
            M_{\frac{  n_j}{\beta J}-\lambda \alpha k}\,f_j\, ,\ M_{\lambda^2\alpha k'}T_{\lambda^2\beta l'} g_0\,
    \rangle
    \notag \\
    &=&\sum_{k,l}\sum_{j=0}^{J-1}
    \langle\,
                   M_{\lambda \alpha k}\,
            T_{\frac{m_j}{\alpha} + \lambda   \beta l J}\,
            P\,
            T_{- \lambda \beta l J  }\,
            M_{\frac{  n_j}{\beta J}-\lambda \alpha k}\,f_j\, ,\ M_{\lambda^2\alpha k'}T_{\lambda^2\beta l'} g_0\,
    \rangle\sigma_{k,lJ+j}
    \notag \\
    &=&  \sum_{k,l}\sum_{j=0}^{J-1} m_{k',l',k,lJ+j}
    \ \sigma_{k,lJ+j} \, .\notag
\end{eqnarray}

In order to use Lemma~\ref{lem:matrixinverse} to show that $M$, and,
therefore, $C\circ \Phi_{\bm f} \circ E$ is not stable, we have to
obtain bounds on the matrix entries of $M$.
Lemma~\ref{lem:elementarybound}, part {\it b}, together with the
rapidly decaying function
\begin{eqnarray}
\widetilde{d}_1=C_{\bm f} \sum_{j=0}^{J-1} T_{\frac{m_j}{\alpha} -
\lambda
                    \beta j} d_1,  \notag 
\end{eqnarray}
 will provide us with
these bounds. In fact, for $k,l,k',l'\in \Z $, we have
\begin{eqnarray}
     |m_{k',l',k´,l''}|&=&|m_{k',l',k´,lJ+j}|\notag \\
       &=&
         \left|
          \langle\,
                            M_{\lambda \alpha k}\,
            T_{\frac{m_j}{\alpha} + \lambda   \beta l J}\,
            P\,
            T_{-\lambda \beta l J  }\,
            M_{\frac{  n_j}{\beta J}-\lambda \alpha k}\,f_j\, ,\
                M_{\lambda^2\alpha k'}T_{\lambda^2\beta l'} g_0\,
            \rangle
         \right|
        \notag \\
        &\leq&
          \langle\,
            T_{  \lambda   \beta (l J+j)}
            \big( T_{\frac{m_j}{\alpha} - \lambda
                    \beta j} \left|P\,
            T_{-\lambda \beta l J  }\,
            M_{\frac{  n_j}{\beta J}-\lambda \alpha k}\,f_j\,\right|\big)\, ,\
               T_{\lambda^2\beta l'} g_0\,
            \rangle
        \notag\\
      &\leq&  \widetilde{d}_1 \ast g_0 \, (\lambda\beta(\lambda l'-l''))\,
            (1+\|(\lambda\beta l J, \frac{  n_j}{\beta J}-\lambda \alpha k ) \|_\infty)^{L_{\bm f}} ,
      \notag
\end{eqnarray}
and
\begin{eqnarray}
     |m_{k',l',k´,l''}|&=&|m_{k',l',k´,lJ+j}|\notag \\
       &=&
         \left|
          \langle\,
                             T_{\lambda \alpha k}\,
            M_{-\frac{m_j}{\alpha} - \lambda   \beta l J}\,
            \big(P\,
            T_{-\lambda \beta l J  }\,
            M_{\frac{n_j}{\beta J}-\lambda \alpha k}\,f_j \,\big)\widehat{\ }\, ,\
                T_{\lambda^2\alpha k'}M_{-\lambda^2\beta l'} g_0\,
            \rangle
         \right|
        \notag \\
  &\leq&
          \langle\,
                 T_{\lambda \alpha k}
            \left|
            \big(P\,
            T_{-\lambda \beta l J  }\,
            M_{\frac{  n_j}{\beta J}-\lambda \alpha k}\,f_j\,\big)\widehat{\ }\right|\, ,\
                T_{\lambda^2\alpha k'} g_0\,
            \rangle
                 \notag \\
  &\leq&  d_2 \ast g_0 (\lambda\alpha(\lambda k'-k))
    (1+\|(\lambda\beta l J, \frac{  n_j}{\beta J}-\lambda \alpha k ) \|_\infty)^{L_{\bm f}} \notag.
\end{eqnarray}
In these calculations, we used that $g_0 \geq0$,
$\widehat{g_0}=g_0$, and
 $g_0(-x)=g_0(x)$, and the Parseval--Plancherel identity.
 Since $\widetilde{d}_1$, $d_2$, and $g_0$
decay rapidly, the same holds for $\widetilde{d}_1\ast g_0$ and
$d_2\ast g_0$. We set
$$
  w(x)= \max\big\{ \widetilde{d}_1 \ast g_0(\lambda\beta x),\
                 \widetilde{d}_1 \ast g_0(-\lambda\beta x),\
                 d_2\ast g_0(\lambda\alpha x),\
                 d_2\ast g_0(-\lambda\alpha x)
          \big\},
$$
and choose a polynomial $p$ of degree $L_{\bm f}$ which satisfies
$$
 (1+\|(\lambda\beta l J, \frac{  n_j}{\beta
J}-\lambda \alpha k ) \|)^{L_{\bm f}}\leq p(\|(k,l)\|_\infty),\quad
j=1,\ldots, J,
$$
 and
obtain $
  |m_{k',l',k,l}|\leq w\big(\max\{|\lambda k'-k|,\ |\lambda l'-l|\}\big)\,p(\|(k,l)\|_\infty)
$
with $w=o\left(x^{-n}\right)$ for $n\in\N$.
Lemma~\ref{lem:matrixinverse} implies that $M$ is not stable, and
therefore $C\circ \Phi_{\bm f} \circ E$ and thus $\Phi_{\bm f}$ are
not stable. \hfill $\square$


\section{Appendix}

{\it Proof of Lemma~\ref{lem:matrixinverse}}

Without loss of generality, we may assume $p(x)= (1+x)^{L}$. First,
we show that if $w:\R_0^+\rightarrow\R_0^+$ with
$w(x)=o\left(x^{-(L{+}2)}\right)$ is monotonically decreasing, then
\eqa{
    K_1^{2L}\sum_{K\geq K_1}K \sum_{k\geq K} k\, w(k)^{2}\to 0\,
    \text{ as } K_1\to \infty . \label{equation:DoubleSumW}
}
This limit 
is proven using the Riemann integral criterium for sums. To this
end, we pick $v\in C_0(\R^+)$ with $w(x)\leq v(x)\,x^{-(L{+}2)}$ and
observe that
    \eqa{  \sum_{K\geq K_1+2}K \sum_{k\geq K} k\, w(k)^{2}
            &\leq&\sum_{K\geq K_1+1}K \sum_{k\geq K+1} k\,
            w(k)^{2}\notag \\
            &\leq&
                \int_{K_1}^\infty x \int_x^\infty y\, w(y)^{2}\,dy \,dx
                \notag \\
            &\leq&  \int_{K_1}^\infty x \int_x^\infty y\, v(y)^{2}
                        y^{-2L-4}\,dy \,dx\notag \\
            &\leq&  \int_{K_1}^\infty x \int_x^\infty  v(y)^{2}
                        y^{-2L-3}\,dy \,dx\notag \\
            &\leq&
                    \frac{\|v|_{[K_1,\infty)}\|_\infty^{2}}{2L+2}
                        \int_{K_1}^\infty x\, x^{-2L-2}\,dx
                     \notag\\
            &\leq&
                    \frac{\|v|_{[K_1,\infty)}\|_\infty^{2}}{2L+2}
                        \int_{K_1}^\infty  x^{-2L-1}\,dx
                    \notag\\
  &\leq&
                \frac{\|v|_{[K_1,\infty)}\|_\infty^{2}}{2L(2L+2)}
  K_1^{-2L}
  =o(K_1^{-2L}).\notag
}
Since $\|v|_{[K_1,\infty)}\|_\infty\to 0$ as $K_1\to \infty$,
(\ref{equation:DoubleSumW}) follows.

Now, we shall use (\ref{equation:DoubleSumW}) to show that $
\inf_{x\in l_0(\Z^2)}\{\frac{\|Mx\|_{l^2}}{\|x\|_{l^2}}\}=0$. To
this end, fix $\epsilon>0$ and note that (\ref{equation:DoubleSumW})
provides us with a $K_1>K_0$ satisfying
$$
  (K_1+3)^{2L}\sum_{K\geq K_1 }
     K  \left( \sum_{k\geq K} k\, w(k)^{2}\right)
\leq  2^{-6}\left(\frac{\lambda -
1}{\lambda}\right)^{2L}\epsilon^{2}\,.
$$
Set $N=\left\lceil \frac{ \lambda (K_1+1)}{\lambda - 1}
\right\rceil$ and $\widetilde{N}=\lceil \frac{N}{\lambda}
\rceil+K_1.$ Then $
  N\leq\frac{\lambda(K_1+2)}{\lambda - 1}
$, and  $
  N\geq\frac{\lambda(K_1+1)}{\lambda - 1}
$ implies $ \lambda N\geq \lambda K_1 + \lambda +N$ and
$$ N\geq
 K_1 + \frac N \lambda +1 > K_1 + \left\lceil \frac N \lambda\right\rceil=\widetilde N.$$
 Therefore, $(2\widetilde{N}+1)^2<(2N+1)^2$ and the matrix
$$
\widetilde{M} = ( m_{j',j})_{\|j'\|_\infty\leq
\widetilde{N},\|j\|\leq N}: \C^{(2N+1)^2}\rightarrow
\C^{(2\widetilde{N}+1)^2}$$ has a nontrivial kernel. We can
therefore choose $\widetilde{x}\in \C^{(2N+1)}$ with
$\|\widetilde{x}\|_{2}=1$ and $\widetilde{M}\widetilde{x}=0$. Define
$x\in l_0(\Z^{2})$ according to $x_j=\widetilde{x}_j$ if
$\|j\|_\infty\leq N$ and $x_j=0$ otherwise, so by construction we
have $\|x\|_{l^2}=1$, and $(Mx)_{j'}=0$ for $\|j'\|_\infty\leq
\widetilde{N}$.

To estimate $(Mx)_{j'}$ for $\|j'\|_\infty> \widetilde{N}$, we fix
$K > K_1$ and one of the $2^3 \big(\lceil\frac{N}{\lambda}
\rceil+K\big)$ indices $j'\in\Z^d$ with
$\|j'\|_\infty=\lceil\frac{N}{\lambda} \rceil+K$. We have $\|\lambda
j'\|_\infty \geq N+K \lambda$ and $\|\lambda j'-j\|_\infty \geq
K\lambda\geq K$ for all $j\in\Z^d$  with $\|j\|_\infty\leq N$.
Therefore
\begin{eqnarray*}
|(Mx)_{j'}|^{2}
  &=&    \Big|\sum_{\|j\|_\infty\leq N} m_{j',j}x_j\Big|^{2}\\
  &\leq& \|x\|^{2}_{2}  \sum_{\|j\|_\infty\leq N}\left|m_{j',j}\right|^{2}\\[.3cm]
  &\leq&\sum_{\|j\|_\infty\leq N}w(\|\lambda j'-j\|_\infty)^{2}
  (1+\|j\|_\infty)^{2L}\\&\leq& (N{+}1)^{2L}\sum_{\|j\|_\infty\leq N}w(\|\lambda
  j'-j\|_\infty)^{2}\\
  &\leq&  (N{+}1)^{2L} \sum_{\|j\|_\infty\geq K}w(\|j\|_\infty)^{2}\\[.3cm]
  &=&   (N{+}1)^{2L} 2^{3}\sum_{k\geq K}k\, w(k)^{2}.
\end{eqnarray*}
Finally, we  compute
\eqa{ \|Mx\|_{l^2}^{2}
  &=& \sum_{j'\in\Z^d}|(Mx)_{j'}|^{2} \notag \\
     &=& \sum_{\|j'\|_\infty \geq \lceil \frac{N}{\lambda} \rceil+K_1}|(Mx)_{j'}|^{2}\notag \\
  &=& \ 2^3 \sum_{\|j'\|_\infty \geq \lceil \frac{N}{\lambda} \rceil+K_1}(N{+}1)^{2L}
   \sum_{k\geq \|j'\|_\infty}k\, w(k)^2\notag\\
  &\leq&  2^6 (N{+}1)^{2L}\sum_{K\geq \lceil \frac{N}{\lambda} \rceil+K_1 }
     K  \sum_{k\geq K} k\, w(k)^2\notag \\
    &\leq&  2^6 \left(\frac{ \lambda (K_1+2)}{\lambda - 1}+1\right)^{2L} \sum_{K\geq \lceil \frac{N}{\lambda} \rceil+K_1 }
       K  \sum_{k\geq K} k\, w(k)^2\notag \\
    &\leq&  2^6 \left(\frac{ \lambda}{\lambda - 1}\right)^{2L} (K_1+3)^{2L} \sum_{K\geq \lceil \frac{N}{\lambda} \rceil+K_1 }
       K  \sum_{k\geq K} k\, w(k)^2
    \leq  \epsilon^{2} \notag
} and obtain $\|Mx\|_{l^2}\leq \epsilon$. Since $\epsilon$ was
chosen arbitrarily and $\|x\|_{l^2}=1$, we have $ \inf_{x\in
l_0(\Z^2)}\{\frac{\|Mx\|_{l^2}}{\|x\|_{l^2}}\}=0$ and $M$ is not
stable. \hfill $\square$

\noindent {\bf Acknowledgement.} The author appreciates discussions
on MIMO channels with Werner Kozek which led to the results in this
paper.

\bibliography{../Bibliography/gabor_goetz}
\bibliographystyle{plain}

Version of \today.
\end{document}